# DIMENSION REDUCTION FOR NONELLIPTICALLY DISTRIBUTED PREDICTORS

By Bing Li[1] and Yuexiao Dong

*Pennsylvania State University*

Sufficient dimension reduction methods often require stringent conditions on the joint distribution of the predictor, or, when such conditions are not satisfied, rely on marginal transformation or reweighting to fulfill them approximately. For example, a typical dimension reduction method would require the predictor to have elliptical or even multivariate normal distribution. In this paper, we reformulate the commonly used dimension reduction methods, via the notion of "central solution space," so as to circumvent the requirements of such strong assumptions, while at the same time preserve the desirable properties of the classical methods, such as $\sqrt{n}$-consistency and asymptotic normality. Imposing elliptical distributions or even stronger assumptions on predictors is often considered as the necessary tradeoff for overcoming the "curse of dimensionality," but the development of this paper shows that this need not be the case. The new methods will be compared with existing methods by simulation and applied to a data set.

**1. Introduction.** Dimension reduction for regression [Li (1991, 1992), Cook and Weisberg (1991), Cook (1994, 1996)] is aimed at finding a lower dimensional vector of linear combinations of the predictors, which retains as much as possible the information in the relationship between the response and the original predictors. Let $X$ be a $p$-dimensional random vector representing the predictor, and let $Y$ be a random variable representing the response. If there is a $p \times q$ ($q \leq p$) matrix $\beta$ such that $Y$ and $X$ are independent conditioning on $\beta^T X$ (henceforth written as $Y \perp\!\!\!\perp X | \beta^T X$), then the column space of $\beta$ is called a dimension reduction space. Under very mild conditions, such as given in Cook (1998), Chiaromonte and Cook (2001) and

Received July 2007; revised February 2008.

[1]Supported in part by NSF Grant DMS-07-04621.

AMS 2000 subject classifications. 62H12, 62G08, 62G09.

*Key words and phrases.* Canonical correlation, central solution spaces, kernel inverse regression, inverse regression, sliced inverse regression, parametric inverse regression.

recently further relaxed by Yin, Li and Cook (2008), the intersection of all such spaces is itself a dimension reduction space. In this case we call the intersection the Central Space and denote it by $S_{Y|X}$ [Cook (1994, 1996)]. A basic problem of dimension reduction is to estimate and make statistical inference about $\mathcal{S}_{Y|X}$.

Commonly used dimension reduction methods, such as those based on inverse conditional moments, require rather strong conditions on the joint distribution of $X$. For example, first-moment-based methods such as Sliced Inverse Regression [Li (1991)] and Ordinary Least Squares [Li and Duan (1989)] require the *linear conditional mean* assumption. That is, $E(X|\beta^T X)$ is a linear function of $X$. Second-moment-based methods such as Sliced Average Variance Estimator [Cook and Weisberg (1991)], Principal Hessian Directions [Li (1992)] and Directional Regression [Li and Wang (2007)] require, in addition, the *constant conditional variance* assumption. That is, $\mathrm{Var}(X|\beta^T X)$ is a nonrandom matrix. Since $\beta$ is unknown, these conditions are assumed to hold for all possible $\beta$. If the first condition holds for all $\beta$, then $X$ has an elliptically-contoured distribution [Eaton (1986)], if both conditions hold for all $\beta$, then $X$ has a multivariate normal distribution. Thus, in effect, either elliptically-contoured or multivariate normal distribution has to be assumed when applying these methods.

If the actual predictors do not satisfy these conditions, current practice often relies on transformation—that is, transform the $p$ components of $X$, $(X_1, \ldots, X_p)$, to $(h_1(X_1), \ldots, h_p(X_p))$ by some functions $h_1, \ldots, h_p$, so that the scatter plot matrix of the transformed predictors resembles that of a multivariate normal distribution. While transformation is a pragmatic—and often effective—strategy, it has both theoretical and practical difficulties. Theoretically, such transformations are intrinsically marginal. It targets the marginal distributions of $X_1, \ldots, X_p$, and as such does not guarantee that $E(X|\beta^T X)$ has desired linearity when $\beta^T X$ is not a set of $X_i$'s. Indeed, there can be hidden nonlinearity among the predictors even if their scatter plot matrix looks perfectly linear. On the other hand, marginal transformations may also be excessive: that $E(X|\beta^T X)$ is linear in $X$ does not require every component of $X$ to be linear against every other component. Practically, whether a transformation has succeeded in transforming a set of observed predictors to an elliptical shape often relies entirely on subjective judgement. Moreover, transforming a high dimensional predictor may be tedious or even infeasible. Another way of dealing with nonellipticity is reweighting [Cook and Nachtsheim (1994)]. However, like transformation, it is not focused on that part of the nonlinearity in the predictors that is relevant to dimension reduction. It is also computationally intensive, especially if the dimension $p$ is high.

When the linear conditional mean and/or constant conditional variance assumptions are satisfied, however, the above-mentioned methods share properties that make them uniquely desirable among nonparametric methods.

First, the slicing (or smoothing) involved in these estimators is over the response $Y$, which is always one-dimensional, regardless of the dimension of $X$. It is well known that smoothing over a high dimensional vector space is undesirable, because the data points within a slice (or a region covered by a smoothing kernel) become sparse at an exponential rate as the dimension increases—a phenomenon often referred to as the "curse of dimensionality" [Bellman (1961)]. Second, the size of the slice (or bandwidth of the kernel) for the above methods need not decrease with the sample size for consistency. These properties make the above methods resemble parametric estimators—they are $\sqrt{n}$-consistent regardless of the dimension of $X$ and have simple asymptotic structure—even though the problems they tackle are in fact nonparametric, in the sense that virtually no assumption is imposed on the conditional distribution of $Y|X$.

In this paper, we introduce a method that does not require linear conditional mean or constant conditional variance, while at the same time preserves all the desirable properties described in the foregoing paragraph. Although the basic idea can potentially apply to methods based on first and second inverse conditional moments (such as SIR, the Sliced Average Variance Estimator and Directional Regression), here we will focus on first inverse conditional moments. The new method is akin to inverse regression, but it is adapted in an automated fashion to the nonlinearity in the predictors, and only that part of the nonlinearity relevant for dimension reduction.

In Section 2, we introduce the key idea of Central Solution Space, a construction that circumvents the linear conditional mean assumption. We study its relation with the Inverse Regression Space and the Central Space. In Section 3 we give a general formulation of inverse regression, which accommodates in a simple form five different dimension reduction methods in the literature and in doing so provides a platform on which to generalize them to the nonelliptical situations. This generalization is then carried out in Section 4. Section 5 is devoted to issues involved in implementation, such as parameterization and optimization. The asymptotic distribution of the estimator is developed in Section 6. Sections 7 and 8 are concerned with simulation comparison and application. Finally, the proofs of the asymptotic results are given in the Appendix.

## 2. Central solution space: The principle.
The best way to explain the central idea of this paper is to explain it in comparison with Sliced Inverse Regression. Assume, without loss of generality, that $E(X) = 0$ and $E(Y) = 0$. Let $\Sigma$ be the covariance matrix of $X$, assumed to be positive definite. Suppose the Central Space has dimension $d$, and let $\beta$ be a $p \times d$ matrix whose columns form a basis in $\mathcal{S}_{Y|X}$. Sliced Inverse Regression is based on the following fact. If

(1) $\qquad E(X|\beta^T X)$ is linear in $X$ (linear conditional mean),

then the random vector $\Sigma^{-1}E(X|Y)$ belongs to $\mathcal{S}_{Y|X}$ almost surely. To see this, let $P(\Sigma)$ be the projection on to $\mathcal{S}_{Y|X}$ with respect to the inner product $\langle a, b \rangle = a^T \Sigma b$ (this will be called the $\Sigma$-inner product); that is, $P(\Sigma) = \beta(\beta^T \Sigma \beta)^{-1}\beta^T \Sigma$. Condition (1) implies $E(X|\beta^T X) = P^T(\Sigma)X$. Hence,

$$\begin{aligned}
\Sigma^{-1}E(X|Y) &= \Sigma^{-1}E[E(X|\beta^T X, Y)|Y] = \Sigma^{-1}E[E(X|\beta^T X)|Y] \\
&= \Sigma^{-1}P^T(\Sigma)E(X|Y) = P(\Sigma)\Sigma^{-1}E(X|Y).
\end{aligned} \tag{2}$$

Thus the random vector $\Sigma^{-1}E(X|Y)$ belongs to the range of the projection operator $P(\Sigma)$, which is $\mathcal{S}_{Y|X}$. Consequently, the column space of the matrix

$$\Sigma^{-1}\operatorname{cov}[E(X|Y)]\Sigma^{-1} \tag{3}$$

is a subspace of $\mathcal{S}_{Y|X}$. This column space will be called the Inverse Regression space, and written as $\mathcal{S}_{\mathrm{IR}}$; the matrix (3) will be written as $A_{\mathrm{IR}}$.

At the first sight, linear conditional mean seems crucial in the foregoing argument. However, note that it is the second equality in (2) that reflects the conditional independence $Y \perp\!\!\!\perp X|\beta^T X$, and it requires virtually no condition. The next two equalities in (2), which require linear conditional mean, merely serve to make $\Sigma^{-1}E(X|Y)$ an explicit vector in $\mathcal{S}_{Y|X}$. This leads us to pay special attention to the equation

$$E(X|Y) = E[E(X|\beta^T X)|Y] \qquad \text{a.s.} \tag{4}$$

*That is, the inverse ($L_2$-) regression of $X$ on $Y$ is the same as the double ($L_2$-) regressions of $X$ on $\beta^T X$ and then on $Y$.* Because of the importance of this equation, we will call it the Inverse Regression Equation. Note that if $\beta$ solves this equation, then so does $\beta A$ for any $d \times d$ nonsingular matrix $A$. That is, the above equation is identified only up to the column space of $\beta$.

DEFINITION 2.1.   If $\beta$ is a matrix of $p$ rows that satisfies the inverse regression equation (4), then span($\beta$) is called a solution space of inverse regression equation.

It is easy to see that if $\beta_1$ satisfies (4) and $\beta_2$ is another matrix such that span($\beta_1$) $\subseteq$ span($\beta_2$), then $\beta_2$ also satisfies (4). For maximum dimension reduction we would like to seek $\beta$ of lowest rank. This leads to the notion of *Central Solution Space*.

DEFINITION 2.2.   If the intersection of any two solution spaces of (4) is itself a solution space of (4), then the intersection of all such spaces will be called the Central Solution Space of the inverse regression equation and written as $\mathcal{S}_{\mathrm{CSS}}$.

By construction, if $\eta$ is a matrix of dimension $p \times d_1$ with $d_1$ greater than the dimension of $\mathcal{S}_{\mathrm{CSS}}$, and if it solves equation (4), then span($\eta$) contains $\mathcal{S}_{\mathrm{CSS}}$.

Central Solution Space is defined under the premise that the intersection of two solution spaces of (4) is again a solution space of (4). The similar premise also underlies the construction of the Central Space, which was recently proved under very weak assumptions by Yin, Li and Cook (2008) in that context. The proof in our context is similar and is omitted.

The next proposition reveals the relation among $\mathcal{S}_{\mathrm{CSS}}$, $\mathcal{S}_{\mathrm{IR}}$ and $\mathcal{S}_{Y|X}$, which is the theoretical foundation of our method. We will say that condition (1) holds for a subspace $\mathcal{S}$ of $\mathbb{R}^p$ if it holds for a matrix $\eta$ whose columns form a basis in $\mathcal{S}$. Henceforth, $P_\eta(\Sigma)$ will denote the orthogonal projection on to span($\eta$) with respect to the $\Sigma$-inner product.

THEOREM 2.1. *Suppose that $Y$ and the elements of $X$ are square integrable and $E(X) = 0$. Then:*

1. $\mathcal{S}_{\mathrm{CSS}} \subseteq \mathcal{S}_{Y|X}$;
2. *If, in addition, condition (1) holds for both $\mathcal{S}_{\mathrm{CSS}}$ and $\mathcal{S}_{\mathrm{IR}}$, then $\mathcal{S}_{\mathrm{IR}} = \mathcal{S}_{\mathrm{CSS}}$.*

PROOF. 1. Let $\beta$ be $p$-row matrix such that span($\beta$) $= \mathcal{S}_{Y|X}$. Then, $Y \perp\!\!\!\perp X | \beta^T X$, which, by (2), implies (4). Thus, $\mathcal{S}_{Y|X}$ is a solution space of (4), and assertion 1 follows.

2. Let $\eta$ be a $p$-row matrix whose columns form a basis in $\mathcal{S}_{\mathrm{CSS}}$. If condition (1) holds for $\eta$, then

$$E(X|Y) = E[E(X|\eta^T X)|Y] = P_\eta^T(\Sigma)E(X|Y) = \Sigma P_\eta(\Sigma)\Sigma^{-1}E(X|Y).$$

Hence, $\Sigma^{-1}E(X|Y) = P_\eta(\Sigma)\Sigma^{-1}E(X|Y)$, and consequently

$$(5) \qquad \Sigma^{-1}\operatorname{Var}[E(X|Y)]\Sigma^{-1} = P_\eta(\Sigma)\Sigma^{-1}\operatorname{Var}[E(X|Y)]\Sigma^{-1}P_\eta^T(\Sigma).$$

Thus, we have $\mathcal{S}_{\mathrm{IR}} \subseteq \mathcal{S}_{\mathrm{CSS}}$.

Conversely, let $\xi$ be a $p$-row matrix whose columns form a basis in $\mathcal{S}_{\mathrm{IR}}$. Then,

$$E\|\Sigma^{-1}E(X|Y) - P_\xi(\Sigma)\Sigma^{-1}E(X|Y)\|^2$$
$$= \operatorname{tr}(A_{\mathrm{IR}}) - \operatorname{tr}[A_{\mathrm{IR}}P_\xi^T(\Sigma)] - \operatorname{tr}[P_\xi(\Sigma)A_{\mathrm{IR}}] + \operatorname{tr}[P_\xi(\Sigma)A_{\mathrm{IR}}P_\xi^T(\Sigma)],$$

where $A_{\mathrm{IR}}$ is as defined in (3). Because span($A_{\mathrm{IR}}$) $=$ span($\xi$) and because $A_{\mathrm{IR}}$ is symmetric, the last three terms on the right (without sign) all reduce to $\operatorname{tr}(A_{\mathrm{IR}})$. Consequently, the above quantity is 0, implying

$$\Sigma^{-1}E(X|Y) = P_\xi(\Sigma)\Sigma^{-1}E(X|Y) \qquad \text{a.s.}$$

Because $E(X|\xi^T X)$ is linear in $\xi$, the right-hand side is

$$P_\xi(\Sigma)\Sigma^{-1}E(X|Y) = \Sigma^{-1}P_\xi^T(\Sigma)E(X|Y) = \Sigma^{-1}E[E(X|\xi^T X)|Y].$$

Hence, $\mathcal{S}_{\mathrm{CSS}} \subseteq \mathcal{S}_{\mathrm{IR}}$.                                          □

Observe that part 1 of the theorem holds without any assumption except the existence of moments; the linearity assumption is required only when $\mathcal{S}_{\mathrm{IR}}$ enters the picture. Thus, if we target $\mathcal{S}_{\mathrm{CSS}}$ instead of $\mathcal{S}_{\mathrm{IR}}$, then we can avoid the linearity assumption.

**3. A general formulation of inverse regression.** Several important dimension reduction methods are directly or indirectly related to the fundamental fact that $\mathcal{S}_{\mathrm{IR}} \subseteq \mathcal{S}_{Y|X}$ under condition (1). These include Ordinary Least Squares (OLS) [Li and Duan (1989)], Sliced Inverse Regression (SIR) [Li (1992)], Parametric Inverse Regression (PIR) [Bura and Cook (2001)], Canonical Correlation [Fung et al. (2002)] and Kernel Inverse Regression (KIR) [Zhu and Fang (1996), Ferre and Yao (2005)]. All these methods rely on the condition (1) for their consistency. The original form of PIR of Bura and Cook (2001) was introduced under the assumption that an inverse parametric regression model is true and under that assumption no restriction needs to be imposed on $X$. However, PIR is in fact consistent when the parametric inverse model is not true, and, in this case, condition (1) is needed for its consistency. This fact is noted in Fung et al. (2002) in a different context. The goal of this paper is to use the general mechanism of the Central Solution Space to extend these methods so that their consistency does not rely on condition (1). For this purpose, we now give a brief outline of the construction of these estimators and synthesize them into a common form.

In the literature, the following estimators are typically described in terms of the standardized predictor. But for our purpose it is easier to describe them in terms of the original predictor (assuming $EX = 0$). This makes no difference at the population level (though it does make a difference at the sample level, where our experience indicates that it is often better to work with standardized predictor).

The OLS estimator is based on the following matrix:

$$A_{\mathrm{OLS}} = \Sigma^{-1}E(YX)E(YX^T)\Sigma^{-1}.$$

Let $\{J_1,\ldots,J_k\}$ be a (measurable) partition of $\Omega_Y$, the sample space of $Y$, and define the discretized version of $Y$ as

$$\delta(Y) = \sum_{\ell=1}^k \ell I(Y \in J_\ell).$$

The SIR estimator is based on the following matrix:

$$A_{\mathrm{SIR}} = \Sigma^{-1} \operatorname{Var}[E(X|\delta(Y))]\Sigma^{-1}.$$

Let $\psi : \mathbb{R}^+ \to \mathbb{R}^+$ be a probability density function $h > 0$ and $y \in \Omega_Y$. Let

$$(6) \qquad \kappa(y, \tilde{y}) = \psi(h^{-1}|y - \tilde{y}|)/E[\psi(h^{-1}|Y - \tilde{y}|)].$$

Because $h$ will be treated as fixed throughout the theoretical development, we suppress the dependence on $h$ from the notation. Let $\tilde{Y}$ be a random variable having the same distribution as $Y$ with $\tilde{Y} \perp\!\!\!\perp (X, Y)$. The KIR estimator is based on the following matrix:

$$(7) \qquad A_{\mathrm{KIR}} = \Sigma^{-1} E\{E[X\kappa(Y, \tilde{Y})|\tilde{Y}]E[X^T\kappa(Y, \tilde{Y})|\tilde{Y}]\}\Sigma^{-1}.$$

Finally, let $h_1, \ldots, h_s$ be square integrable functions from $\Omega_Y$ to $\mathbb{R}$, one of which (say $h_1$) must be taken to be 1 if $Y$ is not centered. Let $H(y) = (h_1(y), \ldots, h_s(y))^T$. Let

$$(8) \qquad \rho(y, \tilde{y}) = H^T(y)E[H(Y)H^T(Y)]^{-1}H(\tilde{y}).$$

The matrix

$$A_{\mathrm{PIR}} = \Sigma^{-1} E\{E[X\rho(Y, \tilde{Y})|\tilde{Y}]E[X^T\rho(Y, \tilde{Y})|\tilde{Y}]\}\Sigma^{-1}$$

is sufficiently general to accommodate (the population versions of) both PIR and Canonical Correlation estimator, though their original forms were quite different. We also note that both estimators allow $Y$ to be a vector, but this is not considered in this paper.

It turns that out all four matrices can be written in the same form, which will greatly simplify the subsequent development and provide insights into the relationship among these methods. Henceforth, for two random elements $U$ and $V$, $U \overset{\mathcal{D}}{=} V$ means that they have the same distribution.

THEOREM 3.1. *The matrices $A_{\mathrm{OLS}}$, $A_{\mathrm{SIR}}$, $A_{\mathrm{KIR}}$, $A_{\mathrm{PIR}}$ can be written in the following form:*

$$(9) \qquad \Sigma^{-1} E\{E[Xg(Y, \tilde{Y})|\tilde{Y}]E[X^T g(Y, \tilde{Y})|\tilde{Y}]\}\Sigma^{-1},$$

*where $g : \Omega_Y \times \Omega_Y \to \mathbb{R}$, $\tilde{Y} \perp\!\!\!\perp (X, Y)$, and $\tilde{Y} \overset{\mathcal{D}}{=} Y$.*

PROOF. That $A_{\mathrm{KIR}}$ and $A_{\mathrm{PIR}}$ have the form (9) follows from their definitions. Also, if we let $g(y, \tilde{y}) = y$, then

$$E(XY) = E(XY|\tilde{Y}) = E[Xg(Y, \tilde{Y})|\tilde{Y}].$$

Thus, $A_{\mathrm{OLS}}$ conforms to (9).

For $A_{\mathrm{SIR}}$, note that, for any $j \in \{1, \ldots, k\}$,

$$E[X|\delta(Y) = j] = \frac{E[XI(\delta(Y) = j)]}{P(\delta(Y) = j)} = \frac{E[XI(\delta(Y) = \delta(\tilde{Y}))|\delta(\tilde{Y}) = j]}{P(\delta(Y) = \delta(\tilde{Y})|\delta(\tilde{Y}) = j)}.$$

Because $Y \stackrel{\mathcal{D}}{=} \tilde{Y}$, the above equality implies that

$$E[X|\delta(Y)] \stackrel{\mathcal{D}}{=} E[XI(\delta(Y) = \delta(\tilde{Y}))|\delta(\tilde{Y})]/P[\delta(Y) = \delta(\tilde{Y})|\delta(\tilde{Y})].$$

Let $g(Y, \tilde{Y}) = I(\delta(Y) = \delta(\tilde{Y}))/P[\delta(Y) = \delta(\tilde{Y})|\delta(\tilde{Y})]$. Then,

$$(10) \qquad\qquad E[X|\delta(Y)] \stackrel{\mathcal{D}}{=} E[Xg(Y, \tilde{Y})|\delta(\tilde{Y})].$$

In the meantime,

$$\tilde{Y} \perp\!\!\!\perp (X, Y) \Rightarrow (\tilde{Y}, \delta(\tilde{Y})) \perp\!\!\!\perp (X, Y)$$

$$\Rightarrow (X, Y) \perp\!\!\!\perp \tilde{Y}|\delta(\tilde{Y}) \Rightarrow (X, Y) \perp\!\!\!\perp \tilde{Y}|\{\delta(\tilde{Y}), \delta(\tilde{Y})\},$$

which, together with $\delta(\tilde{Y}) \perp\!\!\!\perp \tilde{Y}|\delta(\tilde{Y})$, implies that

$$(X, Y, \delta(\tilde{Y})) \perp\!\!\!\perp \tilde{Y}|\delta(\tilde{Y}).$$

See, for example, Dawid ([1979](#)) and Cook ([1998](#)), Proposition 4.6. Hence, the right-hand side of ([10](#)) reduces to $E[Xg(Y, \tilde{Y})|\tilde{Y}]$, and equality ([10](#)) reduces to

$$E[X|\delta(Y)] \stackrel{\mathcal{D}}{=} E[Xg(Y, \tilde{Y})|\tilde{Y}].$$

Thus, $A_{\mathrm{SIR}}$ also has the form ([9](#)).   $\square$

## 4. Extension to nonlinear predictor cases.

The synthesis of the last section provides us a platform on which to extend the five methods to situations where the linear conditional mean condition ([1](#)) does not hold. We now carry out this extension.

4.1. *Central solution spaces.* While Theorem [2.1](#) lays out the basic principle of Central Solution Space, as we have noticed in Section [3](#), various versions of inverse regressions do not take the exact form $\mathrm{Var}[E(X|Y)]$. We now extend Theorem [2.1](#) to accommodate the various forms of inverse regressions, as synthesized in Section [3](#).

Denote the matrix ([9](#)) by $A_{\mathrm{IR}}(g)$ and its column space by $\mathcal{S}_{\mathrm{IR}}(g)$, where $g$ stands for the function $g(Y, \tilde{Y})$ in Theorem [3.1](#). Consider the following equation:

$$(11) \qquad E[Xg(Y, \tilde{Y})|\tilde{Y}] = E[E(X|\beta^T X)g(Y, \tilde{Y})|\tilde{Y}],$$

where, recall that $\tilde{Y} \stackrel{\mathcal{D}}{=} Y$ and $\tilde{Y} \perp\!\!\!\perp (X, Y)$. Let $\mathcal{S}_{\mathrm{CSS}}(g)$ be the Central Solution Space of this equation.

THEOREM 4.1. *Suppose that* $g : \Omega_Y \times \Omega_Y \to \mathbb{R}$ *is a measurable function such that the elements of* $Xg(Y, \tilde{Y})$ *are square integrable. Suppose* $Y$ *and the elements of* $X$ *are square integrable with* $E(X) = 0$ *and* $E(Y) = 0$. *Then:*

1. $\mathcal{S}_{\text{CSS}}(g) \subseteq \mathcal{S}_{Y|X}$;
2. *If, in addition, condition* (1) *holds for both* $\mathcal{S}_{\text{CSS}}(g)$ *and* $\mathcal{S}_{\text{IR}}(g)$, *then* $\mathcal{S}_{\text{IR}}(g) = \mathcal{S}_{\text{CSS}}(g)$.

PROOF. 1. Let $\beta$ be a matrix such that $\text{span}(\beta) = \mathcal{S}_{Y|X}$. Because $\tilde{Y} \perp\!\!\!\perp (X, Y)$, we have

$$\tilde{Y} \perp\!\!\!\perp (X, Y, \beta^T X) \Rightarrow \tilde{Y} \perp\!\!\!\perp (X, Y) | \beta^T X.$$

The expression on the right-hand side, together with $Y \perp\!\!\!\perp X | \beta^T X$, implies that $X \perp\!\!\!\perp Y \perp\!\!\!\perp \tilde{Y} | \beta^T X$, and hence that $X \perp\!\!\!\perp (Y, \tilde{Y}) | \beta^T X$. It follows that $E(X | Y, \tilde{Y}, \beta^T X) = E(X | \beta^T X)$, and consequently

$$
\begin{aligned}
(12) \qquad E[Xg(Y, \tilde{Y}) | \tilde{Y}] &= E[E(X | \beta^T X, Y, \tilde{Y})g(Y, \tilde{Y}) | \tilde{Y}] \\
&= E[E(X | \beta^T X)g(Y, \tilde{Y}) | \tilde{Y}].
\end{aligned}
$$

Thus, $\mathcal{S}_{Y|X}$ is a solution space of (11), and assertion 1 follows.

2. Let $\eta$ be a matrix such that $\text{span}(\eta) = \mathcal{S}_{\text{CSS}}(g)$. Since (1) holds for $\eta$, we have $E(X | \eta^T X) = P_\eta^T(\Sigma) X$, and so (11) becomes

$$E[Xg(Y, \tilde{Y}) | \tilde{Y}] = P_\eta^T(\Sigma) E[Xg(Y, \tilde{Y}) | \tilde{Y}] = \Sigma P_\eta(\Sigma) \Sigma^{-1} E[Xg(Y, \tilde{Y}) | \tilde{Y}],$$

which implies $A_{\text{IR}}(g) = P_\eta(\Sigma) A_{\text{IR}}(g) P_\eta^T(\Sigma)$. Hence, $\mathcal{S}_{\text{IR}}(g) \subseteq \mathcal{S}_{\text{CSS}}(g)$. Conversely, let $\xi$ be a matrix such that $\text{span}(\xi) = \mathcal{S}_{\text{IR}}(g)$. Then,

$$
\begin{aligned}
&E \| \Sigma^{-1} E[(Xg(Y, \tilde{Y}) | \tilde{Y}] - P_\xi(\Sigma) \Sigma^{-1} E[(Xg(Y, \tilde{Y}) | \tilde{Y}] \|^2 \\
&= \text{tr}[A_{\text{IR}}(g)] - \text{tr}[A_{\text{IR}}(g) P_\xi^T(\Sigma)] \\
&\quad - \text{tr}[P_\xi(\Sigma) A_{\text{IR}}(g)] + \text{tr}[P_\xi(\Sigma) A_{\text{IR}}(g) P_\xi^T(\Sigma)].
\end{aligned}
$$

Because $\text{span}[A_{\text{IR}}(g)] = \text{span}(\xi)$ and because $A_{\text{IR}}(g)$ is symmetric, the last three terms on the right (without sign) all reduce to $\text{tr}[A_{\text{IR}}(g)]$. Consequently, the above quantity is 0, implying

$$\Sigma^{-1} E[(Xg(Y, \tilde{Y}) | \tilde{Y}] = P_\xi(\Sigma) \Sigma^{-1} E[(Xg(Y, \tilde{Y}) | \tilde{Y}] \qquad \text{a.s.}$$

Because $E(X | \xi^T X)$ is linear in $X$, the right-hand side reduces to

$$
\begin{aligned}
P_\xi(\Sigma) \Sigma^{-1} E[(Xg(Y, \tilde{Y}) | \tilde{Y}] &= \Sigma^{-1} P_\xi^T(\Sigma) E[(Xg(Y, \tilde{Y}) | \tilde{Y}] \\
&= \Sigma^{-1} E[E(X | \xi^T X) g(Y, \tilde{Y}) | \tilde{Y}].
\end{aligned}
$$

Hence, $\mathcal{S}_{\text{CSS}}(g) \subseteq \mathcal{S}_{\text{IR}}(g)$. $\quad\square$

Let $\mathcal{S}_{\text{OLS}}$, $\mathcal{S}_{\text{SIR}}$, $\mathcal{S}_{\text{KIR}}$, $\mathcal{S}_{\text{PIR}}$ be the columns spaces of $A_{\text{OLS}}$, $A_{\text{SIR}}$, $A_{\text{KIR}}$, $A_{\text{IR}}$. Let $\mathcal{S}_{\text{CSS–OLS}}$, $\mathcal{S}_{\text{CSS–SIR}}$, $\mathcal{S}_{\text{CSS–KIR}}$, $\mathcal{S}_{\text{CSS–PIR}}$ be the column spaces of $A_{\text{IR}}(g)$ with $g$ taken to be the four $g(Y, \tilde{Y})$ functions described in the proof of Theorem 3.1. The following corollary follows immediately from Theorem 4.1.

COROLLARY 4.1. *Suppose all the moments involved in the definitions of* $\mathcal{S}_{\text{OLS}}, \ldots, \mathcal{S}_{\text{PIR}}$ *and* $\mathcal{S}_{\text{CSS–OLS}}, \ldots, \mathcal{S}_{\text{CSS–PIR}}$ *are finite. Then:*

1. $\mathcal{S}_{\text{CSS–OLS}} \subseteq \mathcal{S}_{Y|X}, \mathcal{S}_{\text{CSS–SIR}} \subseteq \mathcal{S}_{Y|X}, \mathcal{S}_{\text{CSS–KIR}} \subseteq \mathcal{S}_{Y|X}, \mathcal{S}_{\text{CSS–PIR}} \subseteq \mathcal{S}_{Y|X}$;
2. *If (1) holds for* $\mathcal{S}_{\text{OLS}}, \ldots, \mathcal{S}_{\text{PIR}}$ *and* $\mathcal{S}_{\text{CSS–OLS}}, \ldots, \mathcal{S}_{\text{CSS–PIR}}$, *then*

$$\mathcal{S}_{\text{OLS}} = \mathcal{S}_{\text{CSS–OLS}}, \qquad \mathcal{S}_{\text{SIR}} = \mathcal{S}_{\text{CSS–SIR}},$$

$$\mathcal{S}_{\text{KIR}} = \mathcal{S}_{\text{CSS–KIR}}, \qquad \mathcal{S}_{\text{PIR}} = \mathcal{S}_{\text{CSS–PIR}}.$$

Again, note that inclusions in part 1 hold without linearity condition (1). Part 2 says that when condition (1) does hold, using Central Solution Space based methods will not lose information as compared to the inverse regression based methods.

4.2. *Objective functions.* We now introduce a population-level objective function whose minimizer yields the solution to (11) for each given $g$. We will also describe how it can be estimated based on an i.i.d. sample of $(X, Y)$. The next theorem will provide a guiding principle for defining the objective function.

THEOREM 4.2. *Suppose that* $\mathcal{S}_{\text{CSS}}(g)$ *has dimension* $d \leq p$, *and let* $\beta$ *be a* $p \times d$ *matrix whose columns form a basis in* $\mathcal{S}_{\text{CSS}}(g)$. *Let* $f(\eta^T X)$ *be a square-integrable function such that, whenever* $\text{span}(\eta) = \text{span}(\beta)$, $f(\beta^T X) = E(X | \beta^T X)$, *and whenever* $\text{span}(\eta) \neq \text{span}(\beta)$,

$$(13) \qquad P\{E[f(\eta^T X)g(Y, \tilde{Y})|\tilde{Y}] \neq E[f(\beta^T X)g(Y, \tilde{Y})|\tilde{Y}]\} > 0.$$

*Let* $\eta_0 \in \mathbb{R}^{p \times d}$ *be the minimizer of*

$$(14) \qquad L(\eta) = E\|E\{[X - f(\eta^T X)]g(Y, \tilde{Y})|\tilde{Y}\}\|^2$$

*over* $\mathbb{R}^{p \times d}$. *Then,* $\text{span}(\eta_0) = \mathcal{S}_{\text{CSS}}(g)$.

PROOF. If $\text{span}(\eta) = \text{span}(\beta)$, then

$$E[f(\eta^T X)g(Y, \tilde{Y})|\tilde{Y}] = E[E(X | \beta^T X)g(Y, \tilde{Y})|\tilde{Y}] = E(Xg(Y, \tilde{Y})|\tilde{Y}) \qquad \text{a.s.}$$

Hence, $L(\eta) = 0$. If $\text{span}(\eta) \neq \text{span}(\beta)$, then, by assumption (13),

$$E\|E\{[f(\eta^T X) - f(\beta^T X)]g(Y, \tilde{Y})|\tilde{Y}\}\|^2 > 0.$$

In the meantime,

$$L(\eta) = E\|E\{[X - f(\beta^T X)]g(Y, \tilde{Y})|\tilde{Y}\}\|^2$$
$$+ E\|E\{[f(\beta^T X) - f(\eta^T X)]g(Y, \tilde{Y})|\tilde{Y}\}\|^2$$
$$+ 2E(E\{[X - f(\beta^T X)]g(Y, \tilde{Y})|\tilde{Y}\}^T$$
$$\times E\{[f(\beta^T X) - f(\eta^T X)]g(Y, \tilde{Y})|\tilde{Y}\}).$$

Because $\mathrm{span}(\beta) = \mathcal{S}_{\mathrm{CSS}}(g)$, the last term is 0. Therefore,

$$L(\eta) \geq E\|E\{[f(\beta^T X) - f(\eta^T X)]g(Y, \tilde{Y})|\tilde{Y}\}\|^2 > 0.$$

Hence, the minimizer of $L(\eta)$ must satisfy $\mathrm{span}(\eta) = \mathrm{span}(\beta)$. □

Rather than assuming $E(X|\beta^T X)$ to be linear in $\beta^T X$ at the outset, as we do for classical methods such as SIR, here we model $E(X|\beta^T X)$ parametrically. Let $f_1, \ldots, f_k$ be functions from $\mathbb{R}^d$ to $\mathbb{R}$. We will assume that $E(X|\beta^T X)$ lies in the space spanned by $f_1(\beta^T X), \ldots, f_k(\beta^T X)$. That is, each component of $E(X|\beta^T X)$ is a linear combination of $f_1(\beta^T X), \ldots, f_k(\beta^T X)$. Under this assumption, the conditional expectation $E(X|\beta^T X)$ can be expressed explicitly as

$$E(X|\beta^T X) = E[XG^T(\beta^T X)]\{E[G(\beta^T X)G^T(\beta^T X)]\}^{-1}G(\beta^T X),$$

where

$$G(\beta^T X) = (f_1(\beta^T X), \ldots, f_k(\beta^T X))^T.$$

Note that we are not assuming—and we do not need to assume—that $E(X|\eta^T X)$ is a linear function of $f_1(\eta^T X), \ldots, f_k(\eta^T X)$ for every $\eta$ in $\mathbb{R}^{p \times d}$. All we need is that this holds at the true $\beta$. We use the function

$$(15) \qquad E[XG^T(\eta^T X)]\{E[G(\eta^T X)G^T(\eta^T X)]\}^{-1}G(\eta^T X)$$

as the $f(\eta^T X)$ in the definition (14) of the objective function $L(\eta)$.

We now construct the sample estimate $L_n(\eta)$ of $L(\eta)$. Suppose that $(X_1, Y_1)$, $\ldots, (X_n, Y_n)$ are an i.i.d. sample of $(X, Y)$. For a function $r(X, Y)$, let $E_n r(X, Y)$ denote the sample average $n^{-1} \sum_{i=1}^n r(X_i, Y_i)$.

1. Center $Y_1, \ldots, Y_n$ and $X_1, \ldots, X_n$ as

$$\hat{Y}_i = Y_i - E_n(Y), \qquad \hat{X}_i = X_i - E_n(X).$$

2. Select $\{f_1, \ldots, f_k\}$ that we deem sufficiently flexible to describe the conditional mean $E(X|\beta^T X)$. For example, based on our experience, it often suffices to include linear and quadratic functions of $\beta^T X$. In this case, the set $\{f_1, \ldots, f_k\}$ includes the following $d(d+3)/2 + 1$ functions

$$\{1\} \cup \{\eta_i^T X : i = 1, \ldots, d\} \cup \{\eta_j^T X \eta_k^T X : 1 \leq j \leq k \leq d\},$$

where $\eta_1, \ldots, \eta_d$ are columns $\eta$. Let

$$\hat{f}(\eta^T \hat{X}) = E_n[\hat{X}G^T(\eta^T \hat{X})]\{E_n[G(\eta^T \hat{X})G^T(\eta^T \hat{X})]\}^{-1}G(\eta^T \hat{X}).$$

3. If using OLS, define $L_n(\eta)$ as

$$E_n\|(\hat{X} - \hat{f}(\eta^T\hat{X}))\hat{Y}\|^2.$$

If using SIR, define $L_n(\eta)$ as

$$\frac{1}{n}\sum_{\ell=1}^{k} E_n[I(\hat{Y} \in J_\ell)]\|E_n[(\hat{X} - \hat{f}(\eta^T\hat{X}))|\hat{Y} \in J_\ell]\|^2,$$

where

$$E_n[(\hat{X} - \hat{f}(\eta^T\hat{X}))|\hat{Y} \in J_\ell] = E_n[(\hat{X} - \hat{f}(\eta^T\hat{X}))I(\hat{Y} \in J_\ell)]/E_n[I(\hat{Y} \in J_\ell)].$$

If using KIR, PIR or the Canonical Correlation estimator, define $L_n(\eta)$ as

$$n^{-1}\sum_{j=1}^{n}\left\|n^{-1}\sum_{i=1}^{n}\{[\hat{X}_i - \hat{f}(\eta^T\hat{X}_i)]g(\hat{Y}_i,\hat{Y}_j)\}\right\|^2,$$

where $g$ is either the function $\kappa$ defined in (6) or the function $\rho$ defined in (8). Note that, for PIR and the Canonical Correlation estimator, $g(\hat{Y}_i, \hat{Y}_j)$ can be factorized into functions of $\hat{Y}_i$ and $\hat{Y}_j$, and thus the above double sum can be simplified as a single sum. We will come back to this in Section 6.

In the following, we refer to the CSS-based modification of a classical estimator as that estimator preceded by the prefix "CSS." For example CSS–SIR is the CSS-modification of SIR.

That the CSS-based methods do not require linearity condition (1) also implies that they are no longer restricted to continuous predictors, because all we need is that $\{f_1(\beta^T X), \ldots, f_k(\beta^T X)\}$ be sufficiently flexible to describe $E(X|\beta^T X)$ whether or not $X$ is continuous. In fact, the application in Section 8 shows that CSS–PIR handles a binary predictor effectively.

## 5. Parameterization of objective function.
In this section, we discuss special issues that arise in the minimization of $L_n(\eta)$. The number

$$E\|(X - E(X|\beta^T X))g(Y, \tilde{Y})\|^2$$

depends on $\beta$ only through its column space and not its specific form. This raises the question of how to parameterize the column space parsimoniously. The similar problem arises frequently in dimension reduction; for example, it arises also in Xia et al. (2002), Cook and Ni (2005) and Cook (2007). The most parsimonious parameterization is via the Grassmann manifold, whose importance in dimension reduction computation is first noted in Cook (2007).

Here, we use a more elementary but rather intuitive parameterization—we assume that columns of $\beta$ to be a set of $d$ orthonormal vectors and parameterize them by the polar coordinate system. First, we represent the class of all orthogonal matrices, denoted by $\mathbb{O}^{p \times p}$, using the polar coordinate system. Note that any matrix in $\mathbb{O}^{2 \times 2}$ can be represented as

$$(16) \qquad \begin{pmatrix} \cos(\alpha) & -\sin(\alpha) \\ \sin(\alpha) & \cos(\alpha) \end{pmatrix}.$$

For an arbitrary dimension $p$, the space $\mathbb{R}^p$ consists of $\binom{p}{2}$ two-dimensional orthogonal hyperplanes, and an orthogonal matrix should be able to rotate a vector along all of them. Thus, any matrix in $\mathbb{O}^{p \times p}$ is the product $\binom{p}{2}$ orthogonal matrices, each resembling the above matrix. In symbols, for $1 \leq i < j \leq p$, let $\theta_{ij}$ be an angle in $[0, \pi]$ and $B_{ij}(\theta_{ij})$ be the matrix in $\mathbb{O}^{p \times p}$ constructed by replacing the $(i,j) \times (i,j)$ submatrix of the identity matrix $I_p$ with the matrix of the form (16). That is,

$$R_{ij}(\theta_{ij}) = \begin{pmatrix} 1 & \cdots & 0 & \cdots & 0 & \cdots & 0 \\ \vdots & & \vdots & & \vdots & & \vdots \\ 0 & \cdots & \cos(\theta_{ij}) & \cdots & -\sin(\theta_{ij}) & \cdots & 0 \\ \vdots & & \vdots & & \vdots & & \vdots \\ 0 & \cdots & \sin(\theta_{ij}) & \cdots & \cos(\theta_{ij}) & \cdots & 0 \\ \vdots & & \vdots & & \vdots & & \vdots \\ 0 & \cdots & 0 & \cdots & 0 & \cdots & 1 \end{pmatrix} \begin{matrix} \\ \\ i\text{th row} \\ \\ j\text{th row} \\ \\ \end{matrix}.$$

$$\phantom{R_{ij}(\theta_{ij}) = } \quad i\text{th column} \qquad j\text{th column}$$

Then, any orthogonal matrix can be represented as

$$(17) \qquad B = \prod_{1 \leq i < j \leq p} B_{ij}(\theta_{ij}).$$

We use the first $d$ columns of $B$ as the polar parameterization of $\eta$.

In this parameterization $\eta$ depends on $\theta$, but not all of them. Let us see what is the subset of the $\theta_{ij}$'s that appear in $\eta$. To begin, consider the case of $p = 5$ and $d = 2$. Then,

$$(18) \qquad B = (B_{12}B_{13}B_{14}B_{15}B_{23}B_{24}B_{25})(B_{34}B_{35}B_{45}),$$

where the parentheses are added artificially to assist discussion. Note that the first 2 columns of $B_{34}, B_{35}, B_{45}$ are the same as those of $I_p$. Therefore, the first 2 columns of $B_{34}B_{35}B_{45}$ are the same as those of $I_p$. This implies that the first 2 columns of $B$ and $B_{12} \cdots B_{25}$ are the same. In other words, $B_{34}B_{35}B_{45}$ can be ignored without changing the first 2 columns of $B$. In general, $\eta$ depends only on the following $\theta_{ij}$'s:

$$(19) \qquad \{\theta_{ij} : 1 \leq i \leq d, i < j \leq p\} \equiv \theta.$$

There are $pd - d(d+1)/2 \equiv m$ parameters in this set. We write the parameterized $\eta$ as $\eta(\theta)$. That is, $\eta(\theta)$ comprises the first $d$ columns of

$$\prod_{1 \leq i \leq d, 1 \leq j \leq p, i < j} B_{ij}(\theta_{ij}).$$

Using this parameterization, we minimize $L_n(\eta(\theta))$ over $\theta$. Because $\theta_{ij}$ and $\theta_{ij} \pm \pi$ give the same direction, we maximize $\theta$ over the set $[0, \pi]^m$. For the initial value of $\theta$, we recommend using the corresponding classical methods such as OLS, SIR, KIR and PIR, or the Outer Product Gradient estimator (OPG) by Xia et al. (2002). Many softwares are available for minimizing functions like $L_n(\eta(\theta))$. For example, the OPTIM function in R works well for our purpose. All it requires is a subroutine that evaluates the objective function and an initial value of $\theta$. We use $\text{span}\{\eta(\hat{\theta})\}$, where $\hat{\theta}$ is the minimizer of $L_n(\eta(\theta))$ as the estimator of $\mathcal{S}_{\text{CSS}}(g)$.

We should note that the polar coordinate system is not the most parsimonious parameterization, in the sense that $\text{span}(\theta)$ does not uniquely determine $\theta$, even though $\theta$ has much lower dimension of $\eta$. Numerically, this causes no difficulty with an appropriately chosen initial value such as described above. However, this does mean that the objective function has a singular Hessian matrix, and this must be taken into account when we derive the asymptotic distribution of $\hat{\theta}$, as we do in the next section.

It is possible for $L_n(\eta(\theta))$ to possess multiple minimizers, and we do occasionally run into this problem. However, this is mitigated by the judicious choice of an initial value. For example, OPG, which is easy to compute, seems to work very well. Furthermore, our experience indicates that as long as $L_n(\eta(\theta))$ is decreased (from the initial value) the performance of the CSS estimators tends to be enhanced regardless of convergence of the algorithm. Thus, if we use a robust minimization algorithm that guarantees to decrease the objective function at each step, the issue of local minima should not cause serious concern. For example, the Nelder–Mead simplex algorithm [Nelder and Mead (1965)], as implemented in OPTIM mentioned previously, is robust in this sense.

**6. Asymptotic distribution.** We now derive the asymptotic distributions of $\hat{\theta}$. Because of limited space, in this paper we will only tackle the asymptotic analysis of CSS–PIR, which includes CSS–OLS as a special case. To further simplify computation, we only consider the case where $h_1(y), \ldots, h_s(y)$ are monomials of $y$ and $f_1(\eta^T X), \ldots, f_k(\eta^T X)$ are monomials of $\eta_j^T X$. Furthermore, we require that both sets of functions must contain the function that is constantly 1. Under this assumption, there is no need to assume $E(X) = 0$ and $E(Y) = 0$ (i.e., no need to center $X$ as $X - EX$ and $Y$ as $Y - EY$) because, for example,

$$\{1, y, \ldots, y^{s-1}\} \quad \text{and} \quad \{1, y - EY, \ldots, (y - EY)^{s-1}\}$$

span the same functional space. The development of the general case is parallel to this simplified case but will have much more terms in the asymptotic expansion, complicating an otherwise transparent argument. Note that this restriction does not apply to estimation, where centering causes no additional complication.

For bookkeeping, we first give a one-to-one correspondence between the double index in (19) and a single index. Let $J = \{(i,j) : j = i+1, \ldots, p, i = 1, \ldots, d\}$. For each $(i,j) \in J$, let

$$t = t(i,j) = p(i-1) - (i-1)i/2 + (j-i).$$

Conversely, for each $t \in \{1, \ldots, m\}$, let

$$i(t) = \max\{i : p(i-1) - (i-1)i/2 \leq t\},$$

$$j(t) = t - [p(i(t)-1) - (i(t)-1)i(t)/2] + i(t).$$

In this arrangement, as the double index $(i,j)$ runs through $J$ with $j$ changing first, the single index $t$ runs through 1 to $m$ and vice versa. Let

$$\phi_{t(i,j)} = \theta_{ij}, \qquad D_{t(i,j)} = B_{ij}.$$

Let $\phi = (\phi_1, \ldots, \phi_m)^T$. The $\eta$ can be equivalently parameterized by $\phi$ as

$$\eta(\phi) = \prod_{t=1}^{m} D_t(\phi_t).$$

Denote the range of $\phi$, $[0, \pi]^m$ by $\Omega_\phi$.

Let $\mathcal{F}$ be a convex class of distributions of $(X, Y)$, which contains the true distribution $F_0$ and all empirical distributions. Let $E_F(\cdot)$ denote the expectation under $F$ and $E(\cdot)$ denote the expectation under $F_0$. We can reexpress $L(\eta(\phi))$ and $L_n(\eta(\phi))$ defined in Section 4.2 as evaluations of a mapping from $\Omega_\phi \times \mathcal{F}$ to $\mathbb{R}$, evaluated at the true distribution $F_0$ and the empirical distribution $F_n$, respectively. Let

$$(20) \quad \begin{aligned} \ell(\phi, F) = \operatorname{tr}\{ & E_F[(X - f(\eta^T(\phi)X))H^T(Y)] \\ & \times [E_F(H(Y)H^T(Y))]^{-1}E_F[H(Y)(X - f(\eta^T(\phi)X))^T]\}. \end{aligned}$$

In this notation, $L(\eta(\phi))$ and $L_n(\eta(\phi))$ becomes $\ell(\phi, F_0)$ and $\ell(\phi, F_n)$.

As we have noted in Section 5, $\phi$ is not uniquely determined by the subspace $\operatorname{span}(\eta(\phi))$. Our asymptotic result reflects this fact by allowing the Hessian matrix of $\ell(\phi, F_0)$ to be singular. Let

$$g(\phi_0, F) = \left[ \frac{\partial \ell(\phi, F)}{\partial \phi} \right]_{\phi = \phi_0}, \qquad W = W(\phi_0, F_0) = \left[ \frac{\partial^2 \ell(\phi, F_0)}{\partial \phi \, \partial \phi^T} \right]_{\phi = \phi_0}.$$

Let $P_W$ be the projection on to the column space of $W$, and let $Q_W = I_m - P_W$. By Taylor expansion, it is easy to see that

$$\ell(\phi_0 + n^{-1/2}\delta, F_0) = n^{-1}\delta^T W\delta + o(n^{-1}),$$

$$\ell(\phi_0 + n^{-1/2}P_W\delta, F_0) = n^{-1}\delta^T W\delta + o(n^{-1}).$$

That is, in a contiguity neighborhood of $\phi_0$, $\ell(\cdot, F_0)$ is unaffected by the component $Q_W\delta$ of the parameter. In other words, locally at $\phi_0$, it is $P_W\delta$ that parameterizes the subspace $\mathrm{span}(\eta(\phi_0 + n^{-1/2}P_W\delta))$, and the component $Q_W\delta$ has no effect on this subspace. Similarly, at the sample level, it can be shown that (not presented here)

$$\ell(\hat{\phi}, F_n) = \ell(\phi_0 + P_W(\hat{\phi} - \phi_0), F_n) + o_p(n^{-1}), \qquad \ell(\hat{\phi}, F_n) = O_p(n^{-1}).$$

Thus, $Q_W(\hat{\phi} - \phi_0)$ has no effect on the sample objective function $\ell(\cdot, F_n)$. For this reason, the asymptotic distribution of relevance is that of $\sqrt{n}P_W(\hat{\phi} - \phi_0)$, rather than that of the full parameter $\sqrt{n}(\hat{\phi} - \phi_0)$. The next theorem gives the asymptotic expansion of $\sqrt{n}P_W(\hat{\phi} - \phi_0)$. Its proof will be given in the Appendix.

THEOREM 6.1.  *Suppose that the regularity conditions described in Section A.1 are satisfied. Let $W^{\dagger}$ be the Moore–Penrose inverse of $W$, and let $g^*(X, Y, \phi_0, F_0)$ be the influence function of the mapping $F \mapsto g(\phi_0, F)$ evaluated at $F_0$. Then,*

$$(21) \qquad P_W(\hat{\phi} - \phi_0) = -W^{\dagger}E_n g^*(X, Y, \phi_0, F_0) + o_p(n^{-1/2}).$$

*The explicit expression for $W$ is given by (32) through (35), and that for $g^*(X, Y, \phi_0, F_0)$ is given by (36) through (38) in the Appendix.*

From expansion (21), we can easily derive the asymptotic distributions of $\sqrt{n}P_W(\hat{\phi} - \phi_0)$.

COROLLARY 6.1.  *Under regularity conditions described in Section A.1,*

$$\sqrt{n}P_W(\hat{\phi} - \phi_0) \xrightarrow{\mathcal{D}} N(0, \Lambda(\phi_0, F_0)),$$

*where $\Lambda(\phi_0, F_0) = W^{\dagger}E\{g^*(X, Y, \phi_0, F_0)[g^*(X, Y, \phi_0, F_0)]^T\}W^{\dagger}$.*

In practice, we can estimate $\Lambda(\phi_0, F_0)$ by replacing $W$ with its sample estimate $W(\hat{\phi}, F_n)$ and replacing $E\{g^*(X, Y, \phi_0, F_0)[g^*(X, Y, \phi_0, F_0)]^T\}$ with

$$n^{-1}\sum_{i=1}^{n}\{g^*(X_i, Y_i, \hat{\phi}, F_n)[g^*(X_i, Y_i, \hat{\phi}, F_n)]^T\}.$$

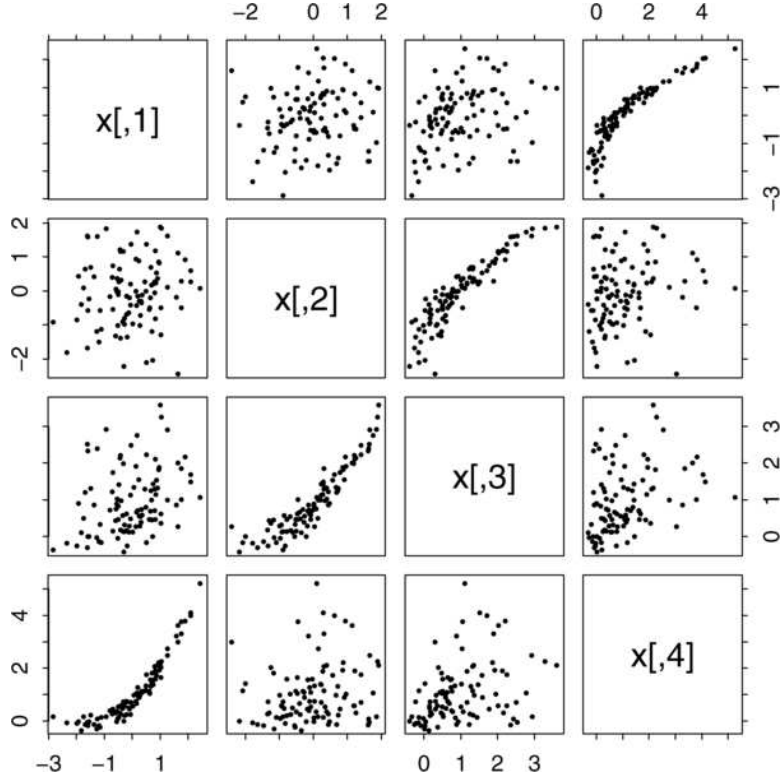

Fig. 1. *Scatter plot matrix for the 4-dimensional, nonelliptically-distributed predictor $X$.*

**7. Simulation comparisons.** In this section, we compare the CSS-based methods with their classical counterparts as well as two adaptive estimators when the predictor $X$ has a nonelliptical distribution. We consider the following three models:

$$\text{Model I:} \quad Y = e^{X_3} + (X_4 + 1.5)^2 + \varepsilon,$$

$$\text{Model II:} \quad Y = 0.4X_3^2 + 3\sin(X_4/4) + 0.5\varepsilon,$$

$$\text{Model III:} \quad Y = X_3/[0.5 + (X_4 + 1.5)^2] + 0.1\varepsilon,$$

where $\varepsilon \sim N(0,1)$ and $\varepsilon \perp\!\!\!\perp X$. We first take the sample size to be $n = 100$. The dimensions of $X$ are chosen to be $p = 4, 6, 8$. Note that, in all three models, $d = 2$ and $\mathcal{S}_{Y|X}$ is spanned by $(0,0,1,\ldots,0)^T$ and $(0,0,0,1,\ldots,0)^T$.

We introduce nonlinearity in the predictor as follows: $X_1 \sim N(0,1)$, $X_2 \sim N(0,1)$,

$$
\begin{aligned}
(22) \qquad & X_3 = 0.2X_1 + 0.2(X_2 + 2)^2 + 0.2\delta, \\
& X_4 = 0.1 + 0.1(X_1 + X_2) + 0.3(X_1 + 1.5)^2 + 0.2\delta,
\end{aligned}
$$

where $\delta \perp\!\!\!\perp (X, Y)$ and $\delta \sim N(0, 1)$. When $p = 6, 8$, $X_5$ through $X_8$ are taken to be independent of $N(0, 1)$ and to be independent of $(X_1, \ldots, X_4)$. Figure 1 shows the scatter plot matrix of $X_1, \ldots, X_4$. Predictors of this type are very common in practice.

We apply three methods based on Central Solution Space CSS–SIR, CSS–PIR and CSS–KIR, as well as their classical counterparts, SIR, PIR and KIR, to the three models. Because CSS–OLS and OLS can only estimate one-dimensional Central Spaces ($d = 1$), we do not include them in the comparison. We also compare with OPG and the Minimum Averaged Variance Estimator (MAVE) introduced by Xia et al. (2002). The simulation sample size is $N = 200$. For SIR and CSS–SIR, the number of slices is taken to be 10, with each slice having equal number of observations. For PIR and CSS–PIR, the function $H(Y)$ is

$$H(Y) = (1, Y, Y^2).$$

For all three CSS methods, the function $G(\eta^T X)$ is taken to be

$$G(\eta^T X) = (1, \eta_1^T X, \eta_2^T X, (\eta_1^T X)^2, \ldots, (\eta_2^T X)^3).$$

For the KIR and CSS–KIR, the function $\psi$ in (6) is taken to be the standard normal density, and the bandwidth $h$ in (6) is taken to be 0.4. The kernel function for OPG and MAVE is taken to be the normal density, with standard deviation (kernel width) taken to be 0.7 for $p = 4, 6$ and 0.8 for $p = 8$. These parameters perform reasonably well in several pilot trial runs.

To assess the accuracy of each method, we use the squared multiple correlation coefficient. Specifically, suppose $U$ and $V$ are $d$ dimensional random vectors, and $\Sigma_{UV}$, $\Sigma_U$ and $\Sigma_V$ are the covariance matrix between $U$ and $V$, the covariance matrix of $U$ and the covariance matrix of $V$, respectively. Then the square multiple correlation coefficient is defined by

$$(23) \qquad \rho^2 = \text{tr}[\Sigma_U^{-1/2} \Sigma_{UV} \Sigma_V^{-1} \Sigma_{VU} \Sigma^{-1/2}] = \text{tr}[\Sigma_V^{-1/2} \Sigma_{VU} \Sigma_U^{-1} \Sigma_{UV} \Sigma^{-1/2}].$$

See Hall and Mathiason (1990). The measure takes maximum value $d$ if $U$ and $V$ have a linear relation and takes minimum 0 if the components of $U$ and $V$ are uncorrelated. At the sample level, given an estimator $\hat{\beta}$ of $\beta$, we use the sample version of the above measure based on

$$\{\hat{\beta}^T X_1, \ldots, \hat{\beta}^T X_n\} \quad \text{and} \quad \{\beta^T X_1, \ldots, \beta^T X_n\}.$$

Note that the larger value of this criterion corresponds to a better dimension reduction estimate.

We compute the errors of estimation by the eight methods, for three models and three choices of $p$ and across the 200 simulated samples. The results are presented in Table 1.

Each entry of Table 1 is formatted as $a(b)$, where $a$ is the average of the above criterion across the 200 simulated samples and $b$ is the standard error of the average. From the table we see that the CSS-based methods are substantially more accurate than their classical counterparts across all 9 cases, indicating that there is much to be gained by correcting the bias caused by nonellipticity. OPG and MAVE perform competently under nonellipticity, but on the whole their improvements are not as sharp as the CSS-based methods. In particular, the accuracy of CSS–KIR dominates that of OPG and MAVE in all 9 cases by substantial margins (relative to the standard deviations). CSS–PIR also performs better than OPG and MAVE in most (8 out of 9) cases. The performance of CSS–SIR is somewhat similar to OPG and MAVE. This is partly due to the fact that slicing is somewhat inefficient, because the inter-slice information is not used—an aspect that cannot be improved by the CSS correction. The loss of intra-slice information by SIR is noticed by Cook and Ni (2006), who proposed a method to reduce it.

For larger sample sizes, the performances of all estimators improve, and MAVE and the CSS-based methods become more similar. Table 2 compares CSS–KIR with KIR, OPG and MAVE for $p = 6$ and $n = 200, 300, 400, 500$.

TABLE 1
*Comparison of CSS and classical estimators*

| Model | Method | $p = 4$ | $p = 6$ | $p = 8$ |
|-------|--------|---------|---------|---------|
| I | PIR | 1.366 (0.017) | 1.336 (0.017) | 1.264 (0.015) |
| | CSS–PIR | 1.658 (0.021) | 1.631 (0.017) | 1.393 (0.017) |
| | SIR | 1.112 (0.013) | 1.100 (0.011) | 1.064 (0.007) |
| | CSS–SIR | 1.735 (0.018) | 1.423 (0.020) | 1.293 (0.019) |
| | KIR | 1.701 (0.014) | 1.661 (0.015) | 1.618 (0.015) |
| | CSS–KIR | 1.832 (0.010) | 1.711 (0.014) | 1.637 (0.017) |
| | OPG | 1.581 (0.023) | 1.377 (0.020) | 1.282 (0.016) |
| | MAVE | 1.785 (0.016) | 1.602 (0.018) | 1.382 (0.017) |
| II | PIR | 1.400 (0.015) | 1.346 (0.015) | 1.349 (0.013) |
| | CSS–PIR | 1.755 (0.018) | 1.558 (0.021) | 1.476 (0.021) |
| | SIR | 1.302 (0.022) | 1.256 (0.017) | 1.208 (0.017) |
| | CSS–SIR | 1.789 (0.013) | 1.439 (0.021) | 1.333 (0.021) |
| | KIR | 1.514 (0.018) | 1.468 (0.016) | 1.437 (0.015) |
| | CSS–KIR | 1.794 (0.015) | 1.551 (0.022) | 1.480 (0.020) |
| | OPG | 1.604 (0.023) | 1.406 (0.023) | 1.302 (0.020) |
| | MAVE | 1.622 (0.022) | 1.397 (0.021) | 1.265 (0.018) |
| III | PIR | 1.149 (0.014) | 1.115 (0.011) | 1.065 (0.009) |
| | CSS–PIR | 1.839 (0.014) | 1.694 (0.018) | 1.557 (0.020) |
| | SIR | 1.265 (0.020) | 1.171 (0.014) | 1.116 (0.013) |
| | CSS–SIR | 1.833 (0.008) | 1.552 (0.020) | 1.454 (0.019) |
| | KIR | 1.146 (0.014) | 1.113 (0.011) | 1.063 (0.009) |
| | CSS–KIR | 1.862 (0.013) | 1.705 (0.019) | 1.613 (0.019) |
| | OPG | 1.742 (0.017) | 1.584 (0.022) | 1.453 (0.020) |
| | MAVE | 1.803 (0.016) | 1.584 (0.021) | 1.375 (0.019) |

TABLE 2
*Comparison KIR, CSS–KIR, OPG and MAVE for larger n's*

| Method | $n = 200$ | $n = 300$ | $n = 400$ | $n = 500$ |
|---|---|---|---|---|
| KIR | 1.704 (0.011) | 1.725 (0.010) | 1.797 (0.005) | 1.781 (0.005) |
| CSS–KIR | 1.816 (0.009) | 1.846 (0.005) | 1.854 (0.004) | 1.861 (0.004) |
| OPG | 1.506 (0.023) | 1.614 (0.022) | 1.681 (0.020) | 1.730 (0.021) |
| MAVE | 1.824 (0.014) | 1.885 (0.012) | 1.847 (0.013) | 1.922 (0.009) |

The kernel width (of $X$) for OPG and MAVE are taken to be $0.6, 0.5, 0.4, 0.4$, and the kernel width (of $Y$) for KIR and CSS–KIR are $0.3, 0.2, 0.1, 0.1$. The basis functions in $H(y)$ now include third polynomials, and the basis functions in $G(\eta^T X)$ include fourth polynomials. We see that, while OPG and KIR still trail behind CSS–KIR, MAVE catches up with CSS–KIR at around $n = 400$ and surpasses it at $n = 500$. This is because, as we can see from (22), the dependence of $X_1$ and $X_2$ on $X_3$ and $X_4$ involves the square root function, and as a consequence $E(X|X_3, X_4)$ does not belong to the polynomials of $\eta^T X$. We have also performed simulation comparisons parallel to those presented in Tables 1 and 2 with $Y$ depending on $X_1$ and $X_2$ instead of $X_3$ and $X_4$, in which case $E(X|\eta^T X)$ does belong to the polynomial family. In this comparison the advantage of the CSS-based methods is more striking, and, for larger $n$, MAVE no longer has the mentioned advantage. These results are not presented for the lack of space.

**8. Application.** We consider data collected for Massachusetts four-year colleges in 1995, which are attempted to study how the percentage of freshmen that graduate (Grad) depends on variables measuring quality of incoming students and features of the colleges. The data is provided as an example data set in MINITAB (release 15, data directory STUDNT12). We restricted attention to $n = 46$ colleges and $p = 8$ predictors, which are: the percentage of freshmen that were among the top 25% percent in their graduating high school class (Top25), the median mathematics SAT score (MSAT), the median verbal SAT score (VSAT), the percentage of applicants accepted by the college (Accept), the percentage of accepted applicants who enroll (Enroll), the student-to-faculty ratio (SFRatio), the out-of-state tuition (Tuition) and whether the college is public or private (PubPriv). Since PIR does not apply to binary data, we first compare PIR and CSS–PIR ignoring the PubPriv variable, and then incorporate PubPriv in the CSS–PIR analysis to see how the latter handles a binary predictor.

The scatter-plot matrix in Figure 2 reveals nonlinearity among predictors—for example, in the relations between Top25 and Accept, Accept and Tuition, VSAT and tuition. The upper panels of Figure 3 present the scatter plots of $Y$ (Grad) versus the first predictors obtained from PIR (left panel) and CSS–PIR (right panel).

Since the true model is unknown, we can no longer use criterion (23) to compare the performances of PIR and CSS–PIR. We will use instead a leave-one-out cross validation criterion to compare their performances [see, e.g., Allen (1974), Stone (1974)]. Let $\tilde{\beta}_{-k}$ and $\hat{\beta}_{-k}$ be the estimated $\beta$ by PIR and CSS–PIR when $(X_k, Y_k)$ is deleted from the sample. From Figure 3, we see that the scatter plots are roughly linear. So, for each of the 46 leave-one-out samples, we fit linear models using both the PIR and the CSS–PIR predictors and predict the deleted $Y_k$ by $\tilde{\beta}_{-k}^T X_k$ and $\hat{\beta}_{-k}^T X_k$ using their respective linear models. The sums of squared prediction errors over the 46 samples for PIR and CSS–PIR are, respectively, 6145 and 5203, indicating a respectable improvement by CSS–PIR.

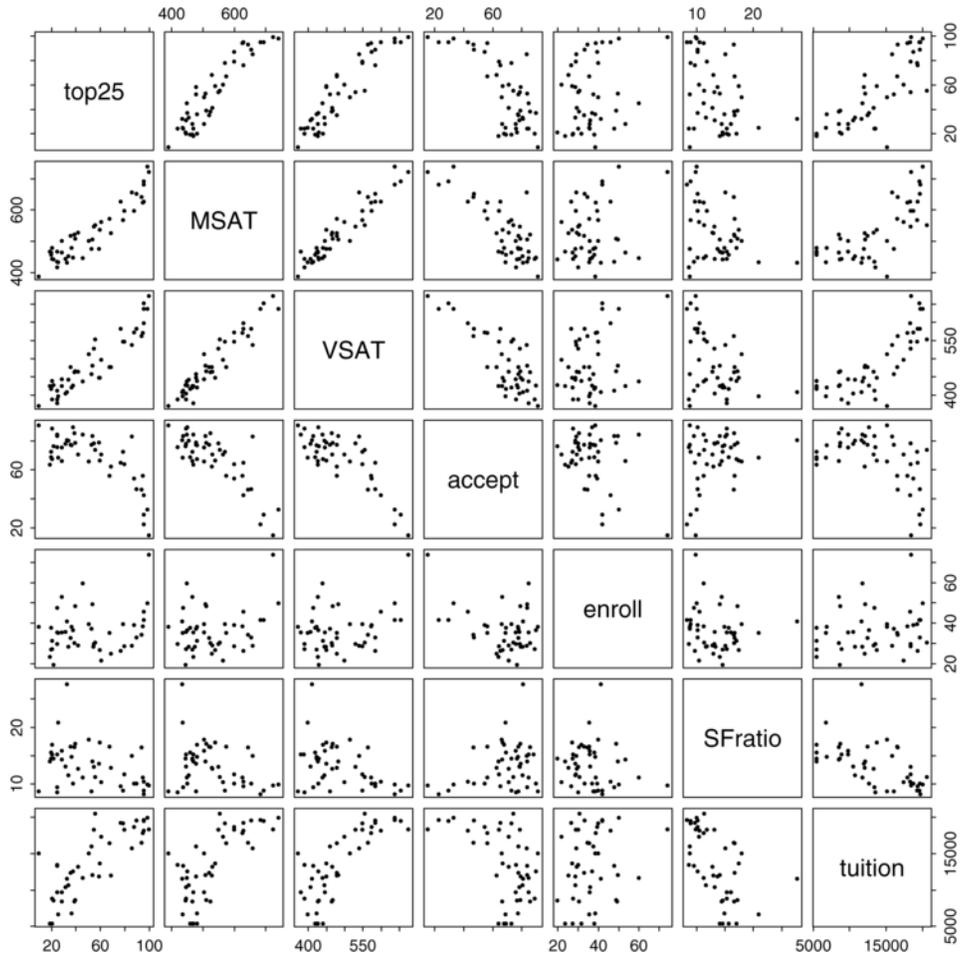

FIG. 2. *Scatter-plot matrix for the seven continuous predictors of the Massachusetts college data.*

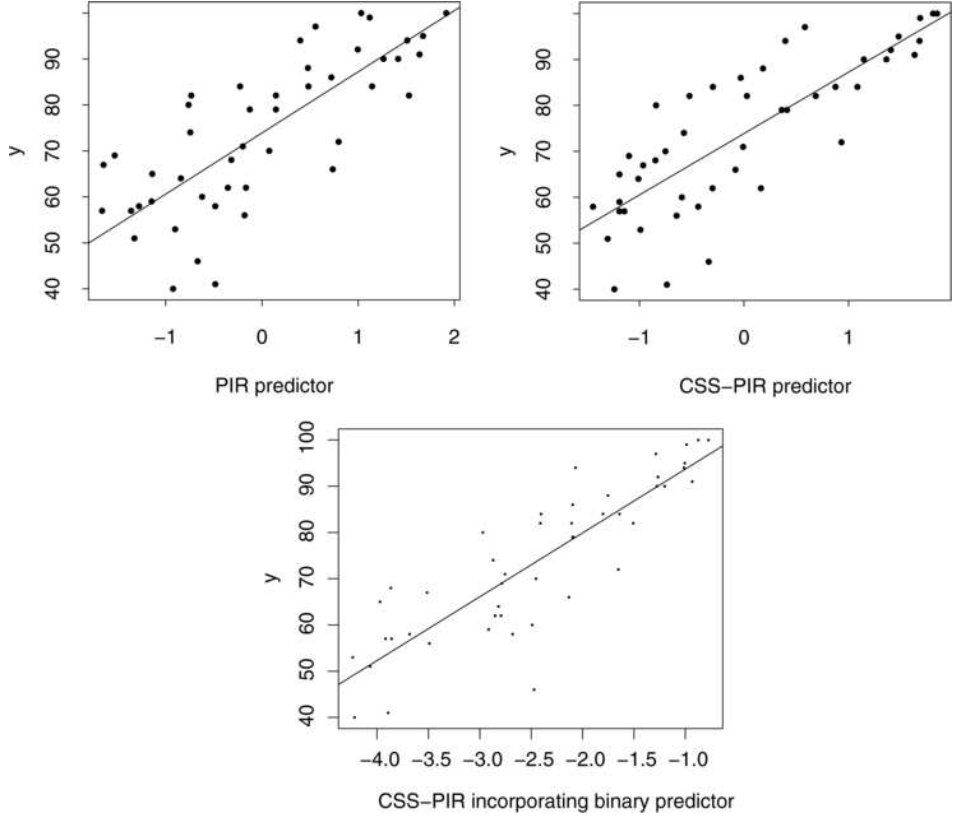

Fig. 3. *Sufficient plots for the Massachusetts college data.*

We now incorporate PubPriv and repeat the CSS–PIR analysis. The 9 public schools are indicated by 1 and the 37 private schools are indicated by 0. The lower panel in Figure 3 is the scatter plot of $Y$ versus the first CSS–PIR predictor after incorporating PubPriv. The cross validation criterion is now further reduced to 4345, another appreciable drop, indicating that CSS–PIR handles the binary predictor effectively.

## APPENDIX: ASYMPTOTIC ANALYSIS

**A.1. Regularity conditions, notation and preliminaries.** The estimator $\hat{\phi}$ in Theorem 6.1 is a function of the empirical distribution $F_n$ of $(X_1, Y_1), \ldots, (X_n, Y_n)$. That is, it has the form $A(F_n)$, where $A$ is a vector. Let $F_0$ be the true distribution of $(X, Y)$ and, for any $\alpha \in [0, 1]$, $F_{n,\alpha} = (1 - \alpha)F_0 + \alpha F_n$.

Then, under regularity conditions,

$$(24) \qquad A(F_n) - A(F_0) = [dA(F_{n,\alpha})/d\alpha]_{\alpha=0} + o_p(n^{1/2}).$$

See von Mises ([1947](#)), Fernholz ([1983](#)) and McCullagh ([1987](#)). When an estimator satisfies ([24](#)), it is called an asymptotically linear estimator [Bickel et al. ([1993](#))]. A wide class of estimators fall into this category. In the following proof, we will assume at the outset that expansion ([24](#)) holds for Theorem [6.1](#). This means that we will omit the proof of consistency and smoothness of the statistical functional $A(\cdot)$. General conditions for estimators defined by minimization of objective functions can be found in van der Vaart ([1998](#)), Chapter 5.

The underlying sufficient condition for expansion ([24](#)) is that the mapping $F \mapsto A(F)$ is Frechet differentiable at $F_0$ (or more generally Hadamard differentiable), with respect to the $\|\cdot\|_\infty$ in a convex family of distributions $\mathcal{F}$ that contains $F_0$ and all empirical distributions. This is not a strong assumption. All estimators discussed in this paper are either themselves functions of sample moments or solutions to equations constructed from sample moments. For example, the key components of SIR and CSS–SIR are the sample conditional moment $E_n[XI(Y \in J_k)]/E_n[I(Y \in J_k)]$, which is a ratio of sample moments. Statistics of this form are typically Frechet differentiable under mild conditions. See, for example, Fernholz ([1983](#)), Chapter 2.

In our context, the Frechet derivative of $A(\cdot)$ at $F_0$ can always be represented as the linear mapping

$$(25) \qquad F \to E_F A^*(X, Y, F_0),$$

where $A^*(X, Y, F_0)$ satisfies $E_{F_0} A^*(X, Y, F_0) = 0$, with its elements belonging to $L_2(F_0)$. When it causes no ambiguity, we will abbreviate $A^*(X, Y, F_0)$ by $A^*(F_0)$. Because of the one-to-one correspondence between the random element $A^*(F_0)$ and mapping ([25](#)), we will refer to $A^*(F_0)$ itself as the Frechet derivative. Moreover, the Frechet derivative, when it exists, coincides with Gateaux derivative, defined as the mapping

$$F \to [D_\alpha A((1 - \alpha)F_0 + \alpha F)]_{\alpha=0}.$$

Hence

$$(26) \qquad [D_\alpha A((1 - \alpha)F_0 + \alpha F)]_{\alpha=0} = E_F A^*(F_0),$$

and consequently the expansion ([24](#)) can be rewritten as

$$(27) \qquad A(F_n) = A(F_0) + E_n A^*(F_0) + o_p(n^{-1/2}).$$

The next lemma provides some basic formulas for Frechet derivatives, defined as the random element $A^*(F_0)$. Let $\rho : \mathcal{F} \to \Theta \subseteq \mathbb{R}^s$ be a Frechet differentiable mapping and $\tau(\rho(F), F)$ be a real, vector or matrix-valued

Frechet differentiable function on $\mathcal{F}$. Let $\rho_0 = \rho(F_0)$. The symbol $\tau^*(\rho_0, F_0)$ denotes the Frechet derivative of the mapping $F \to \tau(\rho_0, F)$ at $F_0$; whereas $\tau^*(F_0, \rho_0(F_0))$ denotes the Frechet derivative of the mapping $F \to \tau(\rho(F), F)$ at $F_0$.

The subsequent expansions demand an efficient notation system for differentiation. We will frequently encounter mappings of the form $Q(\rho(F_\alpha), F_\alpha)$, where $F_\alpha = (1 - \alpha)F_0 + \alpha F$ for some distribution $F \in \mathcal{F}$ and some vector-valued function $\rho(\cdot)$ defined on $\mathcal{F}$. We use $D_\alpha Q(\rho(F_\alpha), F_\alpha)$ to denote the (total) derivative $dQ(\rho(F_\alpha)F_\alpha)/d\alpha$ and $\partial_\alpha Q(\rho, F_\alpha)$ the partial derivative $\partial Q(\rho, F_\alpha)/\partial \alpha$. We use $\partial_{\rho_i} Q(\rho, F_\alpha)$ to denote $\partial Q(\rho, F_\alpha)/\partial \rho_i$, where $\rho_i$ is the $i$ component of $\rho$, and use $\partial_\rho$ to denote the vector of differential operators $(\partial_{\rho_1}, \ldots, \partial_{\rho_m})^T$. We use $\partial_\rho^2$ to denote the matrix $\partial_\rho \partial_\rho^T$. For a single operator such as $\partial_{\rho_i}$, both $\partial_{\rho_i} Q$ and $Q \partial_{\rho_i}$ are to be understood as the derivative $\partial Q/\partial \rho_i$. This is so that differential operation behaves, to a degree, like matrix multiplication. For example, if $q$ is a vector-valued function of $\rho$, then $q \partial_\rho^T$ represents the matrix

$$(q \partial_{\rho_1}, \ldots, q \partial_{\rho_m}) = (\partial_{\rho_1} q, \ldots, \partial_{\rho_m} q)$$

and $\partial_\rho q^T$ denotes the transpose of the above matrix. This notation is helpful in tracking the dimensions of derivative arrays and the ways in which derivatives are arranged in an array. Note, however, that the associative law does not apply: $(q_1 \partial^T)q_2 \neq q_1(\partial^T q_2)$. For this reason, we will always use parentheses to associate a differential array with the function on which it operates.

Lemma A.1. *The following relations hold:*

1. *If the mappings $F \mapsto \rho(F)$ and $F \mapsto \tau(\rho(F), F)$ are Frechet differentiable with respect to $F$ at $F_0$, then*

$$(28) \qquad \tau^*(\rho(F_0), F_0) = \sum_{i=1}^{s} [\partial_{\rho_i} \tau(\rho, F_0)]_{\rho=\rho_0} \rho_i^*(F_0) + \tau^*(\rho_0, F_0);$$

2. *If $\rho(F)$ is a linear functional of $F$, that is, if $\rho(F) = E_F[r(X, Y)]$ for some square-integrable (and vector-valued) function $r$ of $(X, Y)$ that does not depend on $F$, then*

$$\rho^*(F_0) = r(X, Y) - E[r(X, Y)].$$

Proof. Part 2 is well known; [see Fernholz (1983), page 8]. By differentiation,

$$[D_\alpha \tau(\rho(F_\alpha), F_\alpha)]_{\alpha=0} = \sum_{i=1}^{s} [\partial_{\rho_i} \tau(\rho, F_0)]_{\rho=\rho_0} [D_\alpha \rho_i(F_\alpha)]_{\alpha=0} + \partial_\alpha \tau(\rho_0, F_0).$$

By (26), $[D_\alpha \rho_i(F_\alpha)]_{\alpha=0} = E_F \rho_i^*(F_0)$ and $[\partial_\alpha \tau(\rho_0, F_\alpha)]_{\alpha=0} = E_F \tau^*(\rho_0, F_0)$. Hence,

$$[D_\alpha \tau(\rho(F_\alpha), F_\alpha)]_{\alpha=0} = E_F \left[ \sum_{i=1}^{s} [\partial_{\rho_i} \tau(\rho, F_0)]_{\rho=\rho_0} \rho_i^*(F_0) + \tau^*(\rho_0, F_0) \right].$$

By (26) again, the expression inside the brackets on the right-hand side is $\tau^*(\rho(F_0), F_0)$. $\square$

### A.2. Proof of Theorem 6.1. Let

$$
\begin{aligned}
R_1(F) &= E_F[XH^T(Y)], \\
R_2(\phi, F) &= E_F[XG^T(\eta^T(\phi)X)], \\
R_3(\phi, F) &= E_F[G(\eta^T(\phi)X)G^T(\eta^T(\phi)X)], \\
R_4(\phi, F) &= E_F[G(\eta^T(\phi)X)H^T(Y)], \\
R_5(F) &= E_F[H(Y)H^T(Y)].
\end{aligned}
\tag{29}
$$

Note that only $R_2(\phi, F)$, $R_3(\phi, F)$ and $R_4(\phi, F)$ depends on $\phi$. Let

$$R(\phi, F) = R_1(F) - R_2(\phi, F)R_3^{-1}(\phi, F)R_4(\phi, F).$$

Then, by the definition (15) of $f(\eta^T X)$, the $\ell(\phi, F)$ in (20) can be reexpressed as

$$
\begin{aligned}
\ell(\phi, F) &= \mathrm{tr}\{[R_1(F) - R_2(\phi, F)R_3^{-1}(\phi, F)R_4(\phi, F)]R_5^{-1}(F) \\
&\quad \times [R_1(F) - R_2(\phi, F)R_3^{-1}(\phi, F)R_4(\phi, F)]^T\} \\
&= \mathrm{tr}[R(\phi, F)R_5^{-1}(F)R^T(\phi, F)].
\end{aligned}
\tag{30}
$$

Let $\phi(F)$ be the minimizer of $\ell(\phi, F)$. In this notation, $\hat{\phi}$ and $\phi_0$ defined in Section 6 are expressed as $\phi(F_0)$ and $\phi(F_n)$, respectively. Recalling that $g(\phi, F) = \partial_\phi \ell(\phi, F)$, we have

$$g(\phi(F), F) = 0$$

for all $F \in \mathcal{F}$. Take Frechet derivative on both sides, using (28) to obtain

$$W\phi^*(F_0) + g^*(\phi_0, F_0) = 0.$$

Here, $g^*(\phi_0, F_0)$ is to be understood as the Frechet derivative of $F \to g(\phi_0, F)$ at $F_0$ [recall that $g^*(\phi_0, F_0)$ is the abbreviation of $g^*(X, Y, \phi_0, F_0)$]. Multiply both sides of the above equality by $W^\dagger$, and use the fact $W^\dagger W = P_W$ to obtain

$$P_W \phi^*(F_0) = -W^\dagger g^*(\phi_0, F_0), \tag{31}$$

which, by (27), implies (21).

It remains to compute the $m \times m$ nonrandom matrix $W$ and the $m$-dimensional random vector $g^*(\phi_0, F_0)$. Because $R = R_1 - R_2 R_3^{-1} R_4 = 0$ at $(\phi_0, F_0)$, and because $R_5$ does not depend on $\phi$, the $(t, u)$th element of $\partial_\phi^2 \ell(\phi_0, F_0)$ is

$$(32) \qquad W_{tu} = 2 \operatorname{tr}[(\partial_{\phi_t} R) R_5^{-1} (\partial_{\phi_u} R)^T].$$

Here and below, symbols such as $\partial_{\phi_t} R$ and $R_5$ abbreviate functions $\partial_{\phi_t} R(\phi, F)$ and $R_5(F)$ evaluated at $(\phi_0, F_0)$. Because $R_1(F)$ does not depend on $\phi$,

$$(33) \quad \partial_{\phi_t} R = -(\partial_{\phi_t} R_2) R_3^{-1} R_4 + R_2 R_3^{-1} (\partial_{\phi_t} R_3) R_3^{-1} R_4 - R_2 R_3^{-1} (\partial_{\phi_t} R_4).$$

Let $\dot\eta_{\phi_t}$ denote the $p \times d$ derivative matrix $d\eta/d\phi_t$ evaluated at $\phi_0$. Then,

$$\partial_{\phi_t} R_2 = E[(X - EX)(X - EX)^T \dot\eta_{\phi_t} \dot G^T],$$

$$(34) \qquad \partial_{\phi_t} R_3 = E(\dot G \dot\eta_{\phi_t}^T (X - EX) G^T) + E(G(X - EX)^T \dot\eta_{\phi_t} \dot G^T),$$

$$\partial_{\phi_t} R_4 = E(\dot G \dot\eta_{\phi_t}^T (X - EX) H^T).$$

The derivative $\dot\eta_{\phi_t}(\phi)$ can be conveniently computed as follows:

$$(35) \quad \dot\eta_{\phi_t}(\phi) = B_1(\phi_1) \cdots B_{t-1}(\phi_{t-1})[\partial_{\phi_t} B_t(\phi)] B_{t+1}(\phi_{t+1}) \cdots B_m(\phi_m),$$

where $\partial_{\phi_t} B_t(\phi)$ is a $p \times p$ matrix whose $(i(t), j(t)) \times (i(t), j(t))$ submatrix is

$$\begin{pmatrix} -\sin(\phi_t) & -\cos(\phi_t) \\ \cos(\phi_t) & -\sin(\phi_t) \end{pmatrix}$$

and whose other elements are all 0.

We now derive $g^*(\phi_0, F_0)$. Differentiating $\ell(\phi, F)$ with respect to $\phi_t$, and evaluating it at $\phi_0$, we have

$$g_t(\phi_0, F) = 2 \operatorname{tr}\{[\partial_{\phi_t} R(\phi_0, F)] R_5^{-1}(F) R^T(\phi_0, F)\}.$$

Hence, the $t$th component of $g^*(\phi_0, F_0)$ is

$$(36) \qquad g_t^*(\phi_0, F_0) = 2 \operatorname{tr}[(\partial_{\phi_t} R) R_5^{-1} (R^*)^T].$$

By Frechet differentiation of $F \to R(\phi_0, F)$, evaluated at $F_0$, we have

$$(37) \qquad R^* = R_1^* - R_2^* R_3^{-1} R_4 + R_2^* R_3^{-1} R_3^* R_3^{-1} R_4 - R_2 R_3^{-1} R_4^*.$$

Here, symbols such as $R_2^*$ denote the Frechet derivative of the mapping $F \mapsto R_2(\phi_0, F)$ evaluated at $\phi_0$. Using part 2 of Lemma A.1, we deduce that

$$R_1^* = X H^T(Y) - E[X H^T(Y)],$$

$$R_2^* = X G^T(\eta^T(\phi_0) X) - E[X G^T(\eta^T(\phi_0) X)],$$

$$(38) \quad R_3^* = G(\eta^T(\phi_0) X) G^T(\eta^T(\phi_0) X) - E[G(\eta^T(\phi_0) X) G^T(\eta^T(\phi_0) X)],$$

$$R_4^* = G(\eta^T(\phi_0) X) H^T(Y) - E[G(\eta^T(\phi_0) X) H^T(Y)],$$

$$R_5^* = H(Y) H^T(Y) - E[H(Y) H^T(Y)].$$

This completes the proof.

**Acknowledgments.** We thank two referees, an Associate Editor and an Editor for their prompt and extremely helpful reviews which, for example, led us to consider the possibility of using the CSS-based methods to handle discrete predictors.

Department of Statistics
Pennsylvania State University
326 Thomas Building
University Park
Pennsylvania, 16802
USA
E-mail: bing@stat.psu.edu
        yud101@psu.edu